\input amstex

\magnification=1250

\centerline{\bf On Teichm\"uller Space of Surface with Boundary}

\medskip
\centerline{\bf Feng Luo}

\medskip

\centerline{\bf Abstract}
We characterize hyperbolic metrics on compact triangulated surfaces with boundary 
using a variational
principle. As
a consequence, a new parameterization of the  Teichm\"uller space of compact surface with boundary is produced.
In the new parameterization, the Teichm\"uller space becomes
an open convex polytope. It is conjectured that the Weil-Petersson symplectic form can be expressed 
explicitly in terms of the new coordinate.

\medskip

\noindent
\S1. {\bf Introduction}

\noindent
1.1.  The purpose of this paper is to produce a new
parameterization of the Teichm\"uller space of compact surface with non-empty boundary so
 that the lengths of the boundary components are fixed. In this new parameterization, 
the Teichm\"uller space becomes an explicit open convex polytope. Our result can 
be considered as the counter-part of the work of [Le], [Ri] and [Lu1] for  hyperbolic, Euclidean
 and spherical cone metrics on closed surfaces.  In these approaches, constant curvature metrics
 are identified with the critical points of some natural energy functions. 
The energy functions used in [Le] and [Lu1] can be constructed by the cosine laws 
for hyperbolic and spherical triangles. The cosine law for right-angled hyperbolic hexagons 
produces the energy for the current work. All these energies are related to the dilogarithm function. 

As a convention in this paper, all surfaces are assumed to be compact and connected
with non-empty boundary and have negative Euler characteristic unless mentioned otherwise. 
A hyperbolic metric on a compact surface is assumed to have totally geodesic boundary.

\medskip
\noindent
1.2. We begin with a brief recall of the Teichm\"uller spaces. 
Suppose $S$ is a compact surface of non-empty boundary and has negative Euler characteristic. 
It is known that there are hyperbolic metrics with totally geodesic boundary on the surface $S$. Two such 
hyperbolic metrics are \it isotopic \rm if there is an isometry isotopic to the identity between them.
The space of all isotopy classes of hyperbolic metrics  on $S$, denoted by $T(S)$, is called the Teichm\"uller space of the surface $S$.  We are interested in the subspace of $T(S)$ with prescribed boundary lengths. To be precise, let
the boundary components of $S$ be $b_1, ..., b_r$. Assign the i-th boundary component $b_i$ a 
positive number $l_i$ and let $l=(l_1, ..., l_r)$. Then the \it bordered Teichm\"uller space \rm T(S, l) is the subset
of $T(S)$ consisting of those isotopy classes of metrics so that the length of $b_i$ in the metrics is $l_i$. 
The  space $T(S, l)$ has been used recently in calculation of the Weil-Petersson volume 
of the moduli spaces of curves in [Mi].
Using a 3-holed sphere decomposition of the surface $S$ and the associated Fenchel-Nielsen coordinate, it is known
(see [IT] or [Bu]) that $T(S, l)$ is  diffeomorphic to $(\bold R \times \bold R_{>0})^N$ for some integer $N$.

One can decompose the surface $S$  into a union of hexagons instead of 3-holed spheres. These decompositions are called \it ideal triangulations \rm of the surface.
They are also called tri-valent ribbon graphs in the dual setting.
The main result of the paper (theorem 1.2) gives a natural parameterization of the bordered Teichm\"uller space $T(S, l)$ using an ideal triangulation.

\medskip
\noindent
1.3. We now set up the framework by recalling the ideal triangulations and right-angled hyperbolic hexagons.
 A \it colored hexagon \rm
is a hexagon so that three of its non-pairwise adjacent edges are designated as x-edges, the other three edges are
the y-edges.  Take a finite collection of colored hexagons $X$ and identify all y-edges in pairs by homeomorphisms. The quotient space $S$  is a disjoint union of compact surfaces with an \it ideal triangulation \rm  $T$. The
\it edges \rm  and \it 2-cells \rm of the triangulation $T$ are the images of  y-edges and hexagons in $X$ under the quotient map.
 The quotient of each x-edge is called an \it x-arc \rm  in $T$. We use $C(S, T)$,
$E=E(S, T)$, $F=F(S, T)$ to denote the sets of all x-arcs, all edges, and all 2-cells in $T$ respectively.
It is easy to see that every compact surface with negative Euler characteristic and non-empty boundary admits an ideal triangulation.

Suppose $H$ is a colored right-angled hyperbolic hexagon with three y-edges $e_1$, $ e_2$, $ e_3$ and x-edges $f_1, f_2, f_3$ so that $f_i$ is the opposite edge of $e_i$. We call $f_i$ the edge \it facing \rm $e_i$ and $f_j$ an edge
\it adjacent \rm to $e_i$ for $j \neq i$. It is well known that the hexagon $H$ is determined up to isometry preserving coloring by the three lengths $l(e_1), l(e_2), l(e_3)$ of the
y-edges. Furthermore, these three lengths $l(e_1), l(e_2), l(e_3)$ can take any assigned positive numbers
(see [Bu] ).  We define the \it E-invariant \rm of the edge $e_i$ to be the number
 $\frac{1}{2}(l(f_j)+l(f_k) - l(f_i))$ where $\{i,j,k\}=\{1,2,3\}$ and $l(f_i)$ is the length of $f_i$. The E-invariants will
 play the pivotal role in the paper and serve as the coordinate for the bordered Teichm\"uller space $T(S, l)$.

There is a natural one-to-one correspondence between an ideally triangulated compact surface with boundary and a
triangulated closed surface. Namely, for a triangulated closed surface $(S', T')$, let $S$ be the compact
surface obtained from $S'$ by removing a small open regular neighborhood of the union of all vertices. Then
the triangulation $T'$ induces an ideal triangulation $T$ of the surface $S$.  Under this correspondence, vertices
of $T'$ correspond to boundary components of $S$ and
edges of $T'$ correspond to edges of $T$. The 2-cells (hexagons) of $T$ correspond to triangles in $T'$.
The x-arcs in $T$ correspond to angles (or corners) in $T'$. The E-invariant
of an x-arc to be introduced below is the counterpart of the edge invariant introduced by G. Leibon [Le] for hyperbolic
metrics on triangulated closed surface. Here Leibon's invariant assigns
an edge the sum of two angles facing the edge subtracting the sum of the four angles adjacent to the edge.

\medskip
\noindent
1.4. Fix an ideal triangulation $T$ of a compact surface $S$. Each hyperbolic metric $d$  on $S$ produces a \it length function \rm $l_d: E \to \bold R_{>0}$ which assigns each edge $e$ in $T$ the length of the shortest geodesic arc homotopic to $e$ relative to the boundary $\partial S$. It is known that two hyperbolic metrics $d, d'$ on $S$ are isotopic if and only if $l_d =l_{d'}$. Furthermore, any function $l: E \to \bold R_{>0}$ can be realized as $l_d$ for some hyperbolic metric $d$ with totally geodesic boundary.  We call $l_d$ the \it length coordinate \rm
of the metric $d$. Thus the length coordinate parameterizes the Teichm\"uller space 
$T(S)$ by $\bold R_{>0}^E$. However, the image of the bordered Teichm\"uller space $T(S, l)$ 
inside $\bold R_{>0}^E$ is complicated.

The E-coordinate of a hyperbolic metric $d$ on an ideally triangulated surface $(S, T)$ is
defined as follows.  The triangulation $T$ is isotopic to a geometric ideal triangulation $T^*$ in d-metric
such that each edge in $T^*$ is a geodesic segment orthogonal to the boundary $\partial S$. In particular, these
edges $e^*$'s decompose the surface $S$ into a union of right-angled hyperbolic hexagons.
Each edge $e^*$ in $T^*$ is adjacent to one or two hyperbolic hexagons 
(the 2-cells in $T^*$). 
 We define the \it E-invariant \rm of the edge $e$, denoted by $z(e)$, to be the sum of the E-invariants of the corresponding
edge $e^*$ in hyperbolic hexagons adjacent to it. The E-coordinate of the metric is the function $z: E \to \bold R$. Our main results are the following.

\medskip
\noindent
{\bf Theorem 1.1.} \it Suppose $(S, T)$ is a compact ideal triangulated surface. Then,
each hyperbolic metric with totally geodesic boundary on the surface $S$ is determined up to isotopy by its E-coordinate. \rm
\medskip

To state the result for bordered Teichm\"uller space, we have to introduce the notion of \it edge cycle \rm in
the ideal triangulation $T$. A finite collection of ordered edges $e_1, ..., e_k$ in $T$ is said to form a \it cycle \rm
if for each index $i$, counted modulo $k$,  $e_i$ and $e_{i+1}$ are adjacent to some 2-cells in $T$. A \it fundamental
edge cycle \rm is an edge cycle so that each edge in $T$ appears at most twice in the cycle. Each boundary
component of the surface $S$ corresponds to a fundamental edge cycle by counting edges adjacent to the boundary component cyclically. We call these \it boundary edge cycles. \rm

\medskip
\noindent
{\bf Theorem 1.2.} \it Suppose $(S, T)$ is a compact ideally triangulated surface with $r$ boundary components
and $l=(l_1, ..., l_r) \in \bold R_{>0}^r$. Let $E$ be the set of all edges in the triangulation $T$.
Then  E-coordinate is a real analytic diffeomorphism from the bordered Teichm\"uller space $T(S, (l_1, ..., l_r))$ to the convex polytope  $ \{ z: E \to \bold R | $ so that (1.1) and (1.2) hold \}.
Here for each fundamental cycle $e_1, ..., e_k$,
$$ \sum_{i=1}^k z(e_i) > 0, \tag 1.1$$
and for the boundary cycle $e_1, ..., e_k$ corresponding to  the j-th boundary component
$$ \sum_{i=1}^k z(e_i) = l_j. \tag 1.2$$ \rm

\medskip

It seems highly likely that the Weil-Petersson symplectic form on the bordered Teichm\"uller space $T(S, l)$ can be expressed explicitly in terms of the
E-coordinate. See [Bo].

\medskip
\noindent
1.5. The strategy of proving theorems 1.1 and 1.2 goes as follows. By a \it length structure \rm on the ideal
triangulated surface $(S, T)
$ we mean a map $x: C(S, T) \to \bold R_{>0}$ assigning each x-arc a positive number. Length structure is
the counter-part of angle structure on closed triangulated surfaces first introduced by Coin de Verdiere [CV1]
and Rivin [Ri].
The E-invariant of
a length structure $x$ is the function $D_x: E \to \bold R$ assigning each edge $e$ the value $\frac{1}{2}( \sum_{ w \in I} x(w) -\sum_{ w' \in II} x(w'))
$ where $I$ is the set of all x-arcs adjacent to $e$ and $II$ is the set of all x-arcs facing $e$.
 Note that each hyperbolic metric $d$
on $(S, T)$ induces a length structure by measuring the lengths of x-arcs in the hyperbolic ideal triangulation $T^*$
isotopic to $T$. The E-coordinate of the metric $d$ is the E-invariant of its length structure.

For each length structure $x$, we define an energy $V(x)$ by using the cosine law for hyperbolic hexagons.
The energy is a strictly concave function of $x$. Given a function $z: E \to \bold R$, the set of all
length structures with $z$ as E-invariant is a bounded convex set (may be empty) $L(S,T,z)$. We prove that
the maximum point of the strictly concave function $V|: L(S, T, z) \to \bold R$ is exactly the length structure
derived from a hyperbolic metric on the surface $S$. Since a strictly concave function on a  convex set has at most one critical point, this establishes theorem 1.1. To prove theorem 1.2, we show that if  edge invariant $z:E \to
\bold R$ satisfies (1.1) and (1.2), then $L(S, T, z) \neq \phi$ and 
the maximum point of $V|:L(S,T,z) \to \bold R$ always exists.  
The necessity of conditions (1.1) and (1.2) can be verified easily.

We remark that there are now different proofs of theorems 1.1 and 1.2 in [Lu3]. The new proof of theorem 1.1 uses
the Legendre transform of the energy function used in this paper.   Theorem 1.2 can be deduced from theorem 1.1 by
analyzing the map sending the length coordinate to the E-coordinate. It is motivated by Thurston's original proof of
the circle packing theorem in [Th]. 

\medskip
\noindent
1.6. The techniques used in the paper are related to and motivated by the beautiful variational principles
developed  by Colin de Verdiere [CV1],
Br\"agger [Br], Rivin [Ri] and Leibon [Le]  for circle packing, singular Euclidean and singular
hyperbolic structures on surfaces.  
In these works,
the energy functions are all related to the 3-dimensional volume. 
(The energy functional used by Colin de Verdiere was discovered by using the Schlaefli formula for tetrahedra
[CV2].) In [Lu3],
we observe that all these energy functions can be constructed using the cosine law and Legendre transform.
Furthermore, the cosine law produces continuous families of energy functions for variational 
framework on surfaces. As a consequence, theorem 1.1 is one case in a continous family of rigidity theorems.
Whether this is related to the quantum phenomena (for instance, quantum Teichm\"uller theory)
is not clear to us.
Another rich source of energy functions for variational principles on triangulated 
surfaces has been discovered recently
by Bobenko and Springborn [BS] using discrete integrable systems.

Parameterization of Teichm\"uller space using metric ribbon graph has been used extensively recently. 
See for instance the solution of Witten conjecture [Ko] and the stability of the homology of the mapping
 class group [Ha]. In the metric
ribbon graph approach, the key lies in the singular flat metrics arising from Jenkins-Strebel differentials.
The approach in this paper can be considered as a counter-part of metric ribbon graph theory 
using hyperbolic metrics instead of flat metrics. Another approach using hyperbolic geometry has been
worked out by Penner [Pe1], [Pe2].

There are many works on constructing coordinates for Teichm\"uller spaces.
In the work of Bonahon [Bo], a parameterization of the Teichm\"uller space of compact surface 
with boundary was produced using the Bonahon-Thurston shearing cocycles. 
In the work of  Penner [Pe1], he introduced the $\lambda$-coordinate for Teichm\"uller space
of ideally triangulated surfaces with boundary.  The $\lambda$-coordinate and the coordinate introduced
in this paper are quite different.  Other related works are
the papers of M. Schlenker [Sc],  Springborn [Sp] and Ushijima [Us]. 
The relationship between the result in this paper and works of Bonahon, Penner and
Schlenker is not clear to us.  This deserves a further study.
A  fascinating question, suggested by a referee, is
whether there is a geometric interpretation of the energy function used in this paper in terms 
of hyperbolic volume of some hyperideal simplices. See [Sc] for more details.

\medskip
\noindent
1.7. The paper is organized as follows. In section 2, we recall the cosine law and establish some of the basic properties. The energy of a right-angled hyperbolic hexagon is introduced and is shown to be a strictly concave function.  In section 3, we prove theorems 1.1 and 1.2.

\medskip
\noindent
1.8. We thank the referees for careful reading of the paper and for their nice suggestions.
\medskip
\noindent
\S2. {\bf  The cosine law of hyperbolic right-angled hexagon}

\medskip
\noindent
We establish some of the basic properties of the cosine law for hyperbolic right-angled hexagons in this section. In particular, the "capacity" of
a right-angled hyperbolic hexagon is defined. Some of the basic properties of 
the capacity function is established. We do not know the geometric meaning of the capacity.

For simplicity, we assume that the indices $i,j,k$ are  pairwise distinct in this section.
\medskip
\noindent
2.1.   Given a colored hyperbolic right-angle hexagon with  y-edge lengths $y_1, y_2, y_3$,  let $x_1, x_2, x_3$ be the lengths of x-edges
 so that $x_i$-th edge is opposite to the $y_i$-th edge.
The cosine law relating the lengths $x_i$'s with $y_j$'s states that,

$$ \cosh (y_i) =\frac{ \cosh  x_i + \cosh x_j \cosh  x_k}{ \sinh   x_j \sinh   x_k },  \tag 2.1$$
 where $\{i,j,k\}=\{1,2,3\}$.

The partial derivatives of $y_i$ as a function of $x=(x_1, x_2, x_3)$ are given by the following lemma.
\medskip
\noindent
{\bf Lemma 2.1.} \it Let $\{i,j,k\}=\{1,2,3\}$. 

(a) (sine law) $\frac{\sinh(x_i)}{\sinh(y_i)} $ is independent of the index $i$. In particular $A_{ijk}=A_{123}$ where
$A_{ijk} = \sinh(  y_i) \sinh x_j \sinh x_k $,

(b) $\frac{\partial y_i}{\partial x_i} = \frac{\sinh(x_i)}{A_{ijk}}=A \sinh(y_i)$ where $A>0$ is independent of indices,

(c) $\frac{\partial y_i}{\partial x_j} =  - \frac{\partial y_i}{\partial x_i} \cosh y_k$.
 \rm

\medskip
The proof is a simple exercise in calculus, see for instance [Lu2].   

Introduce a new variable $t_i=(x_j+x_k - x_i)/2$ for $\{i,j,k\}=\{1,2,3\}$. Then $x_i = t_j+t_k$. 
 The space of all colored hyperbolic right-angled hexagons parameterized by the new coordinate $t=(t_1, t_2, t_3)$
becomes  $H_3 =\{(t_1, t_2, t_3) \in \bold R^3| t_i+t_j > 0\}$.
We consider $y_i=y_i(t)$ as a smooth function defined on $H_3$. 
\medskip
\noindent
{\bf Corollary 2.2.} \it   The length function $y_i=y_i(t)$ on $H_3$ satisfies,

(a) the differential 1-form $w = \sum_{i=1}^3 \ln \cosh(y_i/2) dt_i$ is closed in the open set $H_3$,

(b) the function $\theta(t) = \int_{(0,0,0)}^t w$ is  strictly concave on $H_3$.

 \rm
\medskip
\noindent
{\bf Remark.} The differential 1-form $w$ in (a) has logarithmic singularity at the point (0,0,0). Thus the integral
in (b) is well defined. This can also be seen in proposition 2.3. 
\medskip
\noindent
{\bf Proof.} To show part (a),  it suffices to prove $\partial ( \ln \cosh(y_i/2) )/\partial t_j$ is symmetric in 
$i \neq j$. By lemma 2.1
and $x_i=t_j+t_k$, 
the partial derivative is found to be,  
$$\frac{ \partial}{\partial t_j}  (\ln \cosh(y_i/2))=\frac{1}{2} \tanh(y_i/2) \frac{ \partial y_i}{\partial t_j}
$$
$$=\frac{1}{2}\tanh(y_i/2)
(\frac{\partial y_i}{\partial x_i} + \frac{\partial y_i}{\partial x_k})$$
$$=\frac{1}{2}\tanh(y_i/2) \frac{\partial y_i}{\partial x_i} (1-\cosh(y_j))$$
$$=\frac{1}{2}  \tanh (y_i/2) A \sinh(y_i) ( 1-\cosh(y_j))$$
$$=\frac{A}{2} \frac{\sinh(y_i/2)}{\cosh(y_i/2)}( 2 \sinh(y_i/2) \cosh(y_i/2) )(-2 \sinh^2(y_j/2))$$
$$=-2A  \sinh^2(y_i/2) \sinh^2(y_j/2), \tag 2.2$$
where $\partial y_i/\partial x_i =A \sinh(y_i)$  is in lemma 2.1(b).
The last expression in (2.2) is symmetric in $i,j$. This establishes (a).

To see part (b), due to part (a) and simply connectivity of $H_3$, the function $\theta(t)$ is well defined. To check the convexity, we calculate the
Hessian of $\theta(t)$. The Hessian matrix is $[ \frac{\partial^2 \theta }{\partial t_r \partial t_s}]_{ 3 \times 3} =
[\frac{\partial}{\partial t_r} ( \ln \cosh(y_s/2))]_{3 \times 3}$.
The diagonal entries of the Hessian can be calculated using lemma 2.1 as follows,

$$\frac{\partial }{\partial t_i}( \ln \cosh(y_i/2)) = \frac{1}{2} \tanh(y_i/2) \frac{\partial y_i}{\partial t_i}$$
$$=\frac{1}{2} \tanh(y_i/2) (\frac{\partial y_i}{\partial x_j} + \frac{\partial y_i}{\partial x_k})$$
$$=\frac{1}{2} \tanh(y_i/2)\frac{\partial y_i}{\partial x_i }( -\cosh(y_k) -\cosh(y_j))$$
$$=\frac{1}{2} \frac{\sinh(y_i/2)}{\cosh(y_i/2)} \sinh(y_i) A ( -2\sinh^2(y_k/2)-2\sinh^2(y_j/2) -2)$$
$$= -2A \sinh^2(y_i/2)( \sinh^2 (y_j/2) + \sinh^2(y_k/2) + 1). \tag 2.3$$

By (2.2), (2.3) and $A>0$,  the  matrix 
$-[\partial^2 \theta/\partial t_r \partial t_s]$ is a
diagonally dominated matrix, i.e.,

$$ -\frac{\partial^2 \theta}{\partial t_i \partial t_i} > | \frac{\partial^2 \theta}{\partial t_i \partial t_j} | + |\frac{\partial^2 \theta}{\partial t_i \partial t_k}|. $$
Since a diagonally dominated matrix is positive definite, 
it follows that the Hessian matrix is negative definite and the function $\theta(t)$ is strictly concave. QED

\medskip

The next proposition relates $\theta(t)$ with the dilogarithm function. Let
$\Lambda_1(u) =\int_0^u \ln \cosh(s) ds$ and $\Lambda_2(u) = \int_0^u \ln \sinh(s) ds$. Both are continuous
in $\bold R$ and related to the dilogarithm function and the Lobachevsky function.

\medskip
\noindent
{\bf Proposition 2.3.} \it The function $\theta(t)$ is
 
$$  2\theta(t) = \Lambda_1(t_1+t_2+t_3) + \sum_{i=1}^3 \Lambda_1(t_i) -\Lambda_2(t_1+t_2) -\Lambda_2(t_2+t_3) -\Lambda_2(t_3+t_1). \tag 2.4$$ \rm

\medskip
\noindent
{\bf Proof. } We will verify that the derivatives of the functions on both sides of (2.4) are the same. 
Note that $ 2\frac{\partial \theta(t)}{\partial t_i} = \ln \cosh^2 (y_i/2)$.
By the cosine law (2.1) and the identity $\cosh^2(u/2) = ( \cosh u + 1)/2$, we have
$$ \cosh^2(y_i/2) =  \frac{ \cosh(x_i) + \cosh(x_j)\cosh(x_k) + \sinh(x_j)\sinh(x_k)}{2 \sinh(x_j)\sinh(x_k)}$$
$$=\frac{ \cosh(x_i) + \cosh(x_j + x_k)}{2 \sinh(x_j)\sinh(x_k)}$$
$$ =\frac{ \cosh((x_1+x_2+x_3)/2) \cosh((x_j+x_k-x_i)/2)}{\sinh(x_j) \sinh(x_k)}$$
$$ =\frac{ \cosh(t_1+t_2+t_3)\cosh(t_i)}{\sinh( t_i+t_j)\sinh(t_i+t_k)}.$$
This shows that
$$ 2\frac{\partial \theta(t)}{\partial t_i } = \ln \cosh(t_1+t_2+t_3) + \ln \cosh(t_i) - \ln \sinh(t_i+t_j) -\ln \sinh(t_i+t_k). \tag 2.5$$
Evidently  the right-hand side of (2.5) is the partial derivative of the right-hand side of (2.4) with respect to the variable
$t_i$. Since both functions
vanish at (0,0,0), this proves the proposition. QED

\medskip

Both functions $\Lambda_1(u)$ and $\Lambda_2(u)$ are continuous in $\bold R$. Thus the function 
$\theta(t)$ has a continuous
extension, still denoted by $\theta(t)$,  to the closure of $H_3$ in $\bold R^3$, i.e., $\theta(t)$ is 
well defined on $\overline{ H_3} =\{(t_1, t_2, t_3) \in \bold R^3 | t_i + t_j \geq 0$, for all
$i \neq j$\}. The next result studies the behavior of the function $\theta(t)$ near the boundary of $H_3$ and near infinity. 

\medskip
\noindent
{\bf Proposition 2.4.}  \it The function $\theta(t)$ defined on $\overline{ H_3}=\{t \in \bold R^3 | t_i + t_j \geq 0\}$ is 
non-negative and bounded. Furthermore, 
for any point $a \in \partial H_3$ and any point $p \in H_3$, 
$$\lim_{s \to 0} \frac{d}{ds} ( \theta( (1-s)a + sp)) = \infty.  \tag 2.6$$
\rm

\medskip
\noindent
{\bf Proof.} Let $f(s) =  2\theta( (1-s) a + sp)$ and $t_i= (1-s)a_i+sp_i$. In the following calculation, the indices are
counted modulo 3.  Then by (2.5),
$$ \frac{df(s)}{ds} = \sum_{i=1}^3  2\frac{\partial \theta}{\partial t_i }(p_i-a_i)$$
$$=\ln \cosh(\sum_{i=1}^3 t_i)\sum_{i=1}^3 (p_i-a_i) +
\sum_{i=1}^3 \ln \cosh(t_i)(p_i-a_i) $$$$-\sum_{i=1}^3 \ln (\sinh(t_{i}+t_{i+1})\sinh(t_i+ t_{i-1}))(p_i - a_i)$$
$$= - \sum_{i=1}^3 \ln\sinh(t_i+ t_{i+1})(p_i+p_{i+1} - a_i - a_{i+1}) +A(s)\tag 2.7$$
where $A(s) =\ln \cosh(\sum_{i=1}^3 t_i)\sum_{i=1}^3 (p_i-a_i) $$+\sum_{i=1}^3 \ln \cosh(t_i)(p_i-a_i) $
so that \newline $\lim_{s \to 0} A(s)$ exists in $\bold R$.

To understand  $\lim_{s \to 0} f(s)$, we will discuss three cases according to the location of the boundary point
$a$: (1) only one of $a_i + a_{i+1}$, for i=1,2,3, is zero; (2) exactly two of three numbers $a_i + a_{i+1}$ are zero,
(3) all of $a_i$'s are zero. Note that $\lim_{s \to 0} (t_i + t_j) = a_i + a_j$.

\noindent
{\bf Case 1}, say $a_1+a_2=0$ and $a_2+ a_3, a_3+a_1>0$. Then by (2.7)  and $\lim_{s \to 0}
(t_i + t_j) >0$ for $(i,j) \neq (1,2)$,
$$\frac{df(s)}{ds} = - \ln \sinh(t_1+t_2) ( p_1+p_2-a_1-a_2) + A_1(s)$$
 where $\lim_{s \to 0} A_1(s)$ exists in $\bold R$.
Due to $p_1+ p_2 > 0$, $a_1+a_2=0$ and $\lim_{s \to 0}\ln \sinh(t_1+t_2) \to -\infty$, it follows that (2.6) holds.

\noindent
{\bf Case 2}, say $a_1+a_2=a_2+a_3=0$ and $a_3+a_1>0$. Then by (2.7), 
$$\frac{df(s)}{ds} = - \ln \sinh(t_1+t_2) ( p_1+p_2-a_1-a_2) -\ln \sinh(t_2+t_3) (p_2+p_3 -a_2 -a_3) + A_2(s)$$
where $\lim_{s \to 0} A_2(s)$ exists in $\bold R$.
Due to $p_i+p_j > a_i+a_j =0$ and $\lim_{s \to 0} (t_i+t_j) =0$ for $(i,j)=(1,2), (2,3)$, it follows again that (2.6) holds.

\noindent
{\bf Case 3}, $a_1=a_2=a_3=0$. Then we have
$$ \frac{df(s)}{ds}= - \sum_{i=1}^3 \ln\sinh(t_i+ t_{i+1})(p_i+p_{i+1} - a_i - a_{i+1}) +A_3(s)$$
where $\lim_{s \to 0} A_3(s)$ exists in $\bold R$.
Since $p_i+p_j > a_i + a_j=0$ and $\lim_{s \to 0} (t_i+t_j) =0$ for all $i,j$, thus (2.6) holds again.

To see that the function $ 2\theta(t)$ is bounded in $\overline{ H_3}$, let us consider for each $u>0$ the
 minimum and maximum values
$m(u)$ and $M(u)$ of $2 \theta (t)$ on the triangle  $X_u=\{ (t_1, t_2, t_3) | t_1+t_2+t_3 = u, t_i + t_j \geq 0\}$.
Since the function $ 2 \theta(t)$ is strictly concave in $X_u$ and  $ 2 \theta(t)$
is symmetric in $t_1, t_2, t_3$, its minimum point is 
achieved at the vertices of $X_u$ and its unique maximum point is invariant under the 
permutations of $t_1, t_2, t_3$. Thus
$M(u) =2\theta(u/3, u/3, u/3)$ and $m(u)=2\theta(u,0,0)$, i.e.,
$M(u) = \Lambda_1(u) + 3\Lambda_1(u/3) - 3\Lambda_2(2u/3)$ and $m(u)=
2 \Lambda_1(u) -2 \Lambda_2(u)$. Since
$m(u) >0$, thus $\theta(t) \geq 0$. On the other hand, $M(u)$ is known to be 
bounded in $[0, \infty)$. Thus  $\theta(t)$ is bounded on $\{t | t_i + t_j \geq 0$ where
$i \neq j$\}. QED

\medskip
\noindent
\S3. {\bf Proofs of Theorems 1.1 and 1.2}

\medskip
Suppose $(S, T)$ is an ideally triangulated surface obtained by identifying y-edges of colored 
hexagons $\tilde{P_1}, ..., \tilde{P_n}$ in pairs
by homeomorphisms $\phi_{ij}$. Let $E=\{e_1, ..., $ \newline
$e_m\}$, $F=\{P_1, ..., P_n\}$ and
$C(S, T) =\{ w_1, ..., w_{3n}\}$ be the sets of all edges, 2-cells, and x-arcs in $T$ respectively. Here the quotient of $\tilde{P_i}$ is
$P_i$.  Each 2-cell $P_i$ contains exactly three x-arcs
$w_{i_1}, w_{i_2}, w_{i_3}$ in $\partial S$. 
We say $w_{i_1}, w_{i_2}, w_{i_3}$ bound the 2-cell and $w_{i_j}$ is an x-arc  of $P_i$. An x-arc $w$ is said to \it facing \rm
(respectively \it adjacent to\rm) an edge $e$ if there is a hexagon $\tilde{P}$ and an x-edge $\tilde{w}$ facing (or adjacent to)
an y-edge $\tilde{e}$ in $\tilde{P}$ so that $w$ and $e$ are the quotients of $\tilde{w}$ and $\tilde{e}$.  If $\tilde{e_1}, \tilde{e_2}, \tilde{e_3}$
are  three y-edges in a hexagon $\tilde{P_i}$, we call their quotient 
edges $e_1, e_2, e_3$ the \it edges \rm of the 2-cell $P_i$. Note
that it may occur $e_1=e_2$. 

Recall that a length structure on $(S, T)$ is a function $x: C(S, T) \to \bold R_{>0}$. Geometrically, a length structure is the same as
a realization of each hexagon $\tilde{P_i}$ by a hyperbolic right-angled hexagon 
(by measuring the lengths of the x-edges). There is no
guarantee that the gluing homeomorphism $\phi_{ij}$ identifies two y-edges of $\tilde{P_k}$'s of the same length. Thus a length structure does not
correspond to a metric on the surface. Hyperbolic metrics on the surface 
$S$ are the same as those length structures so that all $\phi_{ij}$'s identify pairs of edges of the same lengths. 
These length structures are said to be induced from  hyperbolic metrics. 
We consider the
Teichm\"uller space $T(S)$ as the subset of the space of all length structures under this identification. The goal of this paper is to
characterize  $T(S)$ as the critical points of a natural energy function. 

Given a length structure $x: C(S, T) \to \bold R_{>0}$, we define its \it t-coordinate \rm  $t=t_x: C(S, T) \to \bold R$ by
$$ t(w) =\frac{1}{2} ( x(w') + x(w'') - x(w)) \tag 3.1$$
where $w, w', w''$ are the x-arcs in the 2-cell containing $w$.  The length structure $x$ can be recovered from its t-coordinate $t$ by
$x(w) = t(w') + t(w'')$.
The E-invariant of a length structure $x$ is $z=z_x: E \to \bold R$ given by
$$ z(e) = t(w) + t(w') \tag 3.2$$
where $w$ and $w'$ are the x-arcs facing the edge $e$. Note that this definition coincides with the definition of E-invariants introduced in section 1.4
when $x$ is induced from a hyperbolic metric. 

The space of all length structures parameterized by their t-coordinate is $L_t(S, T)=\{t=(t_1, ..., t_{3n}) \in \bold R^3|
t_i+t_j > 0$ whenever x-arcs $w_i$ and $w_j$ are inside a 2-cell in $T$\}.  We define the \it energy  $V(t)$ \rm of
a length structure $t \in L_t(S, T)$ to be,
$$ V(t) = \sum_{ \{w_i, w_j, w_k\} \text{bounds a 2-cell}} \theta(t_i, t_j , t_k) \tag 3.3$$
Geometrically, for a length structure corresponding to a collection of hyperbolic right-angled hexagons, its
 energy is the sum of the
values of $\theta$-function at its hexagons.

Given a function $z: E \to \bold R$, let $L_t(S, T, z) $ be the set of all length structures (in t-coordinates) so that its E-invariant is $z$, i.e.,
$L_t(S, T, z) =\{ t \in \bold R^{3n} | t_i+t_j > 0$ when $w_i, w_j$ are in a 2-cell, and $t_i+t_j=z(e)$ when
x-arcs  $w_i$ and $w_j$ are facing $e$\}.

\medskip
\noindent
{\bf Lemma 3.1.} \it If $L_t(S, T, z) \neq \emptyset$, the energy function $V |: L_t(S,T,z) \to \bold R$ is strictly concave so that
the critical points of $V |$ are exactly the length structures induced from hyperbolic metrics. \rm

\medskip
\noindent
{\bf Proof. } The concavity follows from the concavity of $\theta(t)$. To identify the critical points, we use the Lagrangian multiplier
to $V: L(S, T) \to \bold R$ subject to a set of linear constraints $t_i + t_j = z(e)$ when $w_i, w_j$ are facing $e$. At a critical point $q$
of $V|$, there exists a function $h: E \to \bold R$ (the Lagrangian multiplier) so that for all indices $i$,
$$ \frac{\partial V}{\partial t_i }(q) = h(e) \tag 3.4 $$
where the x-arc $w_i$ is facing the edge $e$.  Suppose the x-arc $w_i$ lies in the 2-cell $P_r$ so that $\tilde{e}$ is
 the y-edge of $\tilde{P_r}$ corresponding to $e$.  We realize all hexagons $\tilde{P_i}$ by hyperbolic right-angled hexagons with x-edge lengths given by the length structure
$q$. Then by corollary 2.2, $\partial V/\partial t_i =\ln \cosh(l(\tilde{e})/2)$. Together with (3.4), this shows that the length $l(\tilde{e})$
of $\tilde{e}$ in the hyperbolic hexagon $\tilde{P_r}$ depends only on the quotient edge $e$ in $T$, i.e., the gluing homeomorphism
$\phi_{ij}$ identify pairs of y-edges of the same hyperbolic lengths. Thus the length structure $q$ is induced from a hyperbolic metric on the surface.
Conversely, suppose we have a length structure $q$ induced from a hyperbolic metric. Then by defining the Lagrangian multiplier $h(e)
$ to be $\ln \cosh (l(e)/2)$, we see that (3.4) holds. Since the constraints are linear functions, it follows that the point $q$ is a critical point of
$V|$.  QED

\medskip
\noindent
3.1. The proof of theorem 1.1 is now simple.  Since a strictly concave function on a convex set has at most one critical point, by lemma 3.1, we see theorem 1.1 holds. 

\medskip
\noindent
{\bf Remark.} Another way to prove theorem 1.1 uses the Legendre transform of the function $\theta$ (see [Lu3]).
By definition,  the Legendre transform $\eta (u_1, u_2, u_3)$ of $\theta(t_1, t_2, t_3)$ is
a strictly concave function in variable $u=(u_1, u_2, u_3)$ where $u_i = \ln \cosh(y_i/2)$ so that 
$\partial \eta/\partial u_i = t_i$.
Now for a hyperbolic metric on the triangulated surface $(S, T)$ with lengths of edges $x=(x_1, ..., x_m)$,  let $u=(\ln \cosh(x_1/2), ..., \ln \cosh(x_m/2))$.  Define $W(u) $ to be the sum of the values of $\eta$ at $(u_i, u_j, u_k)$ where the i-th, j-th and k-th edges bound
a hexagon. Then by definition, $W$ is a smooth strictly concave function so that
the gradient of $W$ is the E-coordinate $z(x)$ of the metric.  It is well known that, for
a smooth strictly  concave function $W$ defined in an open convex set in $\bold R^N$, the gradient $\bigtriangledown W$ is
injective. This implies theorem 1.1.  

\medskip
\noindent
3.2. The proof of theorem 1.2 breaks into two parts. In the first part, we show that if $L_t(S, T, z) \neq \emptyset$, the maximum point
of $V|$ exists in $L_t(S, T, z)$. In the second part, we prove that $L_t(S, T, z) \neq \emptyset$ if and only if condition (1.1) in theorem 1.2 holds.

\medskip
\noindent
3.3. To prove the first part, by proposition 2.3, the function $V: L_t(S, T) \to \bold R$ can be extended continuously to the closure $\overline{L_t(S, T)}$
of $L_t(S,T) \subset \bold R^{3n}$.  On the other hand, the set $L_t(S,T,z)$ is bounded. Indeed, we have,

\medskip
\noindent
{\bf Lemma 3.2.} \it Suppose $e_1, ..., e_k$ form an edge cycle in $T$. Then
$$ \sum_{i=1}^k z(e_i) = \sum_{i=1}^k x(w_{n_i}) \tag 3.5$$
where $w_{n_i}$ is the x-arc adjacent to both $e_i$ and $e_{i+1}$ with indices counted modulo $k$.  In particular,
$\sum_{i+1}^k z(e_i) > 0$ for all edge cycles. If the length structure $x$ is induced from a hyperbolic metric, then
$$ \sum_{i=1}^k z(e_i) = l_j \tag 3.6$$
for the boundary cycle $e_1, ..., e_k$ associated to the j-th boundary component of length $l_j$. \rm

\medskip
\noindent
{\bf Proof. } The proof is a simple calculation using the following identity. Namely, the sum of two t-coordinates
$t_i =\frac{1}{2}( x_j + x_k -x_i)$ and $t_j=\frac{1}{2}(x_i+x_k-x_j)$ is $x_k$ where $x_k$ is the edge adjacent to both $y_i$-th and $y_j$-th edge.
Thus (3.5) follows from the above identity and the definition of the edge cycles. The identity (3.6) follows from (3.5) and the definition
of the boundary length. QED

\medskip
\noindent
{\bf Corollary 3.3.} \it (a) The space $L_t(S, T, z)$ is bounded.

(b) If $z : E \to \bold R$ is an E-invariant associated to a hyperbolic metric, then (1.1) and (1.2) hold. \rm

\medskip
\noindent
{\bf Proof.} Part (b) follows from (3.5) and (3.6). To see part (a), we consider the x-coordinate of length structures. Take a length structure $x: C(S, T) \to \bold R_{>0}$ with 
E-invariant $z$. Each x-arc $w$ is in some boundary component $b_i$ of the surface. Thus by (3.5)
$$ 0\leq x(w) \leq \sum_{ w' \subset b_i} x(w') =\sum_{j=1}^k z(e_{n_j}) $$
where $e_{n_1}, ..., e_{n_k}$ is the boundary edge cycle associated to $b_i$. It shows that $x(w)$ is
bounded. QED

\medskip
By corollary 3.3, the energy function $V|$ can be extended continuously to the compact 
closure $\overline{ L_t(S, T,z)}$. In particular, it
has a maximum point $p=(p_1, ..., p_{3n})$ (considered as a t-coordinate) in $\overline{ L_t(S,T,z)}$. We claim that the maximum point $p$ is in $L_t(S, T,z)$. 

We prove the claim by contradiction. Suppose otherwise, by definition, there are pairs of x-arcs, say
 $w_1$ and $w_2$ in a 2-cell so that $p_1+p_2=0$.
A 2-cell $P$ in $T$ is said to be \it degenerated \rm with respect to $p$ if there are two x-arcs $w_i, w_j$ 
 in $P$ so that $p_i+p_j=0$. Let $I$ be
the set of all degenerated 2-cells and $II$ be the set of all non-degenerated 2-cells. Take a point  $q \in L_t(S, T, z)$. For each 2-cell $P$ in 
the triangulation with x-arcs $w_i, w_j ,w_k$, consider the limit
$$h(P)= \lim_{s \to 0} \frac{d}{ds}( \theta((1-s)p_i+sq_i, (1-s)p_j + sq_j, (1-s)p_k+sq_k)).$$
By proposition 2.3, this limit is finite if $P \in II$ and is the positive infinite if $P \in I$.  By the assumption that $I  \neq \emptyset$,  it follows
that
$$ \lim_{s \to 0} \frac{d}{ds}(V((1-s)p+sq))=\sum_{P \in I} h(P) + \sum_{P \in II} h(P) =\infty.  \tag 3.7$$

On the other hand, since $p$ is the maximum point, the function $V((1-s)p + sq)$ has a maximum point at $s=0$. Thus
$\limsup_{s \to 0} \frac{d}{ds}(V((1-s)p + sq)) \leq 0$. This is a contradiction to (3.7).

By the claim, $p \in L_t(S, T, z)$.  By lemma 3.1, it follows that $p$ is induced by a hyperbolic metric.
To summary, we have shown that if $L_t(S, T, z) \neq \emptyset$, then there exists a hyperbolic metric with E-invariant $z$.

\medskip
\noindent
3.4. The necessity of conditions (1.1) and (1.2) follows from corollary 3.3(b). To finish the proof of theorem 1.2, it remains to show the following.

\medskip
\noindent
{\bf Lemma 3.4.} \it Given a function $z : E \to \bold R$ so that (1.1) holds, then $L(S, T, z) \neq \emptyset$. \rm

\medskip
\noindent
{\bf Proof}. Let us consider length structures parameterized by the x-coordinate. Here $x: C(S, T) \to \bold R$ so that $x_i =x(w_i)$ and
$x=(x_1, ..., x_{3n})$. Let the set of all edges be $E =\{e_1, ..., e_m\}$ and $z_i = z(e_i)$.
By definition, $L(S, T, z) =\{ x \in \bold  R^{3n} | $ so that (3.8) and (3.9) hold\},
$$ \sum_{ i \in I} x_i -\sum_{i \in J} x_j = 2 z(e)  \quad  \text{ for each edge $e$} \tag 3.8$$
where $\{w_i | i\in I\}$ and $\{w_j | j \in J\}$ are the sets of x-arcs adjacent to and  facing  the edge $e \in E$ respectively, and
$$ x_i > 0, \quad \text{ for all i.} \tag 3.9$$ 
Consider the set $D=\{ (y_1, ..., y_m) \in \bold R^m |   y_i + y_j \geq y_k$ whenever $e_i, e_j, e_k$ form the edges of a 2-cell\} and the linear programming problem
$\min \{ \sum_{i=1}^n y_i z_i |  $ $y \in D$ \}. The dual linear programming problem is
 $\max\{ 0  | x \in L(S, T, z)\}$ by the construction.
By the duality theorem of linear programming (in fact Farkas lemma suffices in this case) (see [BL] for instance), 
$L(S, T, z) \neq \emptyset$ if and only if for each non-zero vector $y \in D$, $\sum_{i=1}^n y_i z_i  >0$.

The set $D$ is
a cone in $\bold R^m$. Furthermore, the inequalities $y_i + y_j \geq y_k$ and $y_k+y_i \geq y_j$
imply that $y_i \geq 0$. Thus $D =\{ y=(y_1, ..., y_m) \in \bold R^m | y_i \geq 0, $ and $  y_i + y_j \geq y_k$ whenever $e_i, e_j, e_k$ form the edges of a 2-cell\}.
This shows that $D$ can be identified with the space of all measured laminations on the surface $S$ where $y_i$'s are the geometric intersection coordinates. By the work of [Th] (see  [Mo] or others), it is known that every measured lamination on $S$ considered as a vector in $D$
is a non-negative linear combination of those vectors in $D$ associated to essential simple loops.
Furthermore, these essential simple loops can be assumed to  intersect each edge $e \in E$ in at most two points.
It follows that each of these simple loop corresponds to a 
fundamental edge cycle by counting the edges intersecting  it.  In particular, each fundamental cycle $c=(e_{n_i}, ..., e_{n_k})
$  in the triangulation $T$ corresponds to base vector $v_c =(y_1,..., y_n)$ where $y_i=0$ if $i \neq n_j$ and $y_i=1$ if $i=n_j$.
The above discussion shows that each vector in the cone $D$ is a non-negative linear combination of the base vectors $v_c$'s. 
Now condition (1.1) says that $\sum_{i=1}^n y_iz_i$ is positive at every base vector $v_c$. It follows that
for all $y \in D-\{0\}$, $\sum_{i=1}^n y_i z_i  >0$.   QED

\bigskip

\centerline{\bf References}

\bigskip 

[Bo]  Bonahon, F., Shearing hyperbolic surfaces, bending pleated surfaces and Thurston's symplectic form. Ann. Fac. Sci. Toulouse Math. (6) 5 (1996), no. 2, 233--297.

[BL] Borwein, Jonathan M.; Lewis, Adrian S., Convex analysis and nonlinear optimization.
Springer-Verlag, New York, 2000.

[Bu]  Buser, Peter, Geometry and spectra of compact Riemann surfaces. Progress in Mathematics, 106. Birkhauser Boston, Inc., Boston, MA, 1992.

[Br] Br\"agger, W.,  Kreispackungen und Triangulierungen. Enseign. Math., 38:200-217,1992.

[BS]
Bobenko, Alexander I.; Springborn, Boris A., Variational principles for circle patterns and Koebe's theorem. Trans. Amer. Math. Soc. 356 (2004), no. 2, 659--689 

[CV1] Colin de Verdiere, Yves, Un principe variationnel pour les empilements de cercles. Invent. Math. 104 (1991), no. 3, 655--669. 

[CV2] Colin de Verdiere, Yves, private communication.

[Ha] Harer, John L., Stability of the homology of the mapping class groups of orientable surfaces. Ann. of Math. (2) 121 (1985), no. 2, 215--249.

[IT]  Imayoshi, Y.; Taniguchi, M., An introduction to Teichm\"uller spaces. Translated and revised from the Japanese by the authors. Springer-Verlag, Tokyo, 1992.

[Ko] Kontsevich, Maxim, Intersection theory on the moduli space of curves and the matrix Airy function. Comm. Math. Phys. 147 (1992), no. 1, 1--23. 

[Le] Leibon, Gregory, Characterizing the Delaunay decompositions of compact hyperbolic surfaces. Geom. Topol. 6 (2002), 361--391 

[Lu1] Luo, Feng,  A Characterization of spherical polyhedral surfaces, to appear in Jour. Diff. Geom.,
 http://front.math.ucdavis.edu/math.GT/0408112, 2004

[Lu2] Luo, Feng,  Volume and angle structures on 3-manifolds, preprint, 2005,
\newline http: front.math.ucdavis.edu/math.GT/0504049

[Lu3] Luo, Feng,  Rigidity of polyhedral surfaces, in preparation.

[Mi] Mirzakhani, Maryam, Weil-Petersson volumes and intersection
theory on the moduli space of curves. preprint, 2004.

[Mo]   Mosher, Lee, Tiling the projective foliation space of a punctured surface. Trans. Amer. Math. Soc. 306 (1988), no. 1, 1--70. 

[Pe1] Penner, R. C., Decorated Teichm\"uller theory of bordered surfaces.
 Comm. Anal. Geom. 12 (2004), no. 4, 793--820. 

[Pe2]  Penner, R. C., The decorated Teichm\"uller space of punctured surfaces. Comm.
 Math. Phys. 113 (1987), no. 2, 299--339.

[Ri1] Rivin, Igor, Euclidean structures on simplicial surfaces and hyperbolic volume. Ann. of Math. (2) 139 (1994), no. 3, 553--580

[Sc] Schlenker, Jean-Marc, Hyperideal polyhedra in hyperbolic manifolds,
\newline http://front.math.ucdavis.edu/math.GT/0212355.

[Sp] Springborn, B., A variational principle for weighted Delaunay triangulations and hyperideal polyhedra,
http://front.math.ucdavis.edu/math.GT/0603097

[Th] Thurston, William, Geometry and topology of 3-manifolds, lecture notes, Math Dept., Princeton University,   1978,                                       
 at www.msri.org/publications/books/gt3m/

[Us] Ushijima, Akira, A canonical cellular decomposition of the Teichm\"uller space of compact surfaces with boundary. Comm. Math. Phys. 201 (1999), no. 2, 305--326.
\medskip
\noindent

Department of Mathematics

Rutgers University

Piscataway, NJ 08854, USA

\medskip

Center of Mathematical Science

Zhejiang University

Hangzhou, China

email: fluo\@math.rutgers.edu

\end

\medskip
\medskip

For my own record, the linear programming problem used in lemma 3.4 is the following

Lemma. \it $ min\{ y.z | Ay \geq 0\} \geq 0$ iff $\{ x | A^t x=z, x \geq 0\} \neq \emptyset$. \rm

\medskip
Here is the proof. If the set is non-empty, say x is in the set, then $y.z = y. A^t.x = Ay.x \geq 0$.
Conversely, if $min \geq 0$, we prove that the set is non-empty by contradiction. Suppose the set
is empty, then $z $ is not in the image of $\bold R_{\geq 0}^n$ under $A^t$. By the separation theorem,
there is a vector $y$ so that $y . z <0$ and $ y . A^t x \geq 0$ for all $x \geq 0$. This shows that
$Ay . x \geq 0$ for all $ x \geq 0$, i.e., $Ay \geq 0$. We have obtained a contradiction that
$min <0$. QED

\end